 \documentclass[brochure,english,12pt]{bourbaki}

 \newcommand{\R}{\mathbb{R}}
 
 \newcommand{\E}{\mathcal{E}}
 
 \def\mean#1{\mathchoice
         {\mathop{\kern 0.2em\vrule width 0.6em height 0.69678ex depth -0.58065ex
                 \kern -0.8em \intop}\nolimits_{\kern -0.4em#1}}%
         {\mathop{\kern 0.1em\vrule width 0.5em height 0.69678ex depth -0.60387ex
                 \kern -0.6em \intop}\nolimits_{#1}}%
         {\mathop{\kern 0.1em\vrule width 0.5em height 0.69678ex
             depth -0.60387ex
                 \kern -0.6em \intop}\nolimits_{#1}}%
         {\mathop{\kern 0.1em\vrule width 0.5em height 0.69678ex depth -0.60387ex
                 \kern -0.6em \intop}\nolimits_{#1}}}

 \renewcommand{\AA}{\mathscr{A}}

 \newcommand{\HH}{\mathcal{H}}

 \renewcommand{\tt}{{\mbox{\boldmath$t$}}}
 \newcommand{\uu}{{\mbox{\boldmath$u$}}}

 \newcommand{\tauV}{{\kern-3pt\tau}}

 \newcommand{\oVVVk}{\overline{\mbox{\boldmath$V$}}\kern-3pt}
 \newcommand{\tVVVk}{\tilde{\mbox{\boldmath$V$}}\kern-3pt}

 \newcommand{\<}{\langle}
 \renewcommand{\>}{\rangle}
 \renewcommand{\a}{\alpha}
 \renewcommand{\b}{\beta}

 \newcommand{\e}{\varepsilon}

 \renewcommand{\l}{\lambda}

 \newcommand{\n}{\nabla}

 \renewcommand{\O}{\Omega}

 \renewcommand{\i}{\infty}
 \newcommand{\p}{\partial}

 \renewcommand{\det}{\operatorname{det}}

\date{Juin 2018}
\bbkannee{70\textsuperscript{e} ann\'ee, 2017--2018}
\bbknumero{1147}
\title{On the Monge-Amp\`ere equation}
\subtitle{}
\author{Alessio FIGALLI}
\address{ETH Z\"urich, Department of Mathematics\\
R\"amistrasse 101, 8092 Z\"urich, Switzerland}
\email{alessio.figalli@math.ethz.ch}

\begin{document}

\maketitle

\tableofcontents

\section{Introduction}

The Monge-Amp\`ere equation is a nonlinear partial differential equation arising in several problems from analysis and geometry, such as the prescribed Gaussian curvature equation, affine geometry, optimal transportation, etc.

In its classical form, this equation is given by
\begin{equation}
\label{eq:MAclassical}
\det D^2u=f(x,u,\n u) \qquad \text{in $\Omega$},
\end{equation}
where $\Omega\subset \R^n$ is some open set,
$u:\Omega \to \R$ is a convex function,
and $f:\Omega\times \R\times \R^n\to \R^+$ is given.
In other words, the Monge-Amp\`ere equation prescribes the 
product of the eigenvalues of the Hessian of $u$, in contrast with the ``model'' elliptic equation
$\Delta u=f$ which prescribes their sum.
As we shall explain later, the convexity of the solution $u$ is a necessary condition to make the equation degenerate elliptic, and therefore to hope for regularity results.

The goal of this note is to give first a general overview of the classical theory, and then discuss some recent important developments on this beautiful topic.

\section{Historical background}
The Monge-Amp\`ere equation draws its name from its initial formulation in two dimensions by the French mathematicians Monge \cite{Mon2} and Amp\`ere \cite{Am}.

The first notable results on the existence and regularity for the Monge-Amp\`ere equation are due to Minkowski \cite{M1,M2}:
by approximating a general bounded convex set with convex polyhedra with given faces areas, he proved the existence of a weak solution to the ``prescribed Gaussian curvature equation" (now called ``Minkowski problem'').
Later on, using convex polyhedra with given generalized curvatures at the vertices, Alexandrov also proved the existence of a weak solution
in all dimensions, as well as the $C^1$ smoothness of solutions in two dimensions \cite{AL1,AL2,AL3}.

In higher dimension, based on his earlier works, Alexandrov \cite{AL4} (and also Bakelman \cite{Bak1} in two dimensions) introduced
a notion of generalized solution to the Monge-Amp\`ere equation and proved the existence and uniqueness of solutions to the Dirichlet problem (see Section \ref{sect:Alex sol}).
The notion of weak solutions introduced by Alexandrov
(now called ``Alexandrov solutions'') has  continued to be frequently used in recent years,
and a lot of attention has been drawn to prove smoothness of
Alexandrov solutions under suitable assumptions on the right hand side 
and the boundary data.

The regularity of weak solutions in high dimensions is a very delicate problem.
For $n \geq 3$, Pogorelov found a convex function in $\R^n$ which is not of class $C^2$ but satisfies the
Monge-Amp\`ere equation in a neighborhood of the origin with positive analytic right hand side (see \eqref{eq:example pogorelov} below).
It became soon clear that the main issue in the lack of regularity was the presence of 
a line segment in the graph of $u$.
Indeed, Calabi \cite{cal smooth} and Pogorelov \cite{pog smooth} were able to prove a priori 
 interior second and third derivative estimate for strictly convex solutions, or for solutions which do not contain a line segment with both endpoints on boundary. However, in order to perform the computations needed to deduce these a priori estimates, $C^4$ regularity of the solution was needed. Hence, a natural way to prove existence of smooth solutions was to approximate the Dirichlet problem with nicer problems for which $C^4$ solutions exist, apply Pogorelov and Calabi's estimates to get $C^2/C^3$ a priori bounds, and then take the limit in the approximating problems. 
 This argument was successfully implemented by Cheng and Yau \cite{chengyau} and Lions \cite{lions} to obtain the interior smoothness of solutions.

Concerning boundary regularity, thanks to the regularity theory developed by Ivochkina \cite{Ivo},
Krylov
\cite{Kry2}, and
Caffarelli-Nirenberg-Spruck \cite{CNS}, one may use the continuity method 
and Evans-Krylov's estimates \cite{Ev,Kry} to obtain globally smooth solutions to the Dirichlet problem (see Section \ref{sect:continuity}).
In particular,
 Alexandrov  solutions are smooth up to the boundary provided all given data are smooth. 
 
In all the situations mentioned above, one assumes that $f$ is positive and sufficiently smooth. When $f$ is merely 
bounded away from zero and infinity, Caffarelli proved the $C^{1,\alpha}$ regularity of strictly convex solutions \cite{Caf2}.
Furthermore, when $f$ is continuous (resp. $C^{0,\alpha}$), using perturbation arguments
Caffarelli proved interior $W^{2,p}$ estimate for any $p > 1$ (resp. interior $C^{2,\alpha}$ estimates)
\cite{Cafw2p}.

As explained in Section \ref{sect:applic classic}, these results can be applied to obtain both the regularity in the Minkowski problem and in the optimal transportation problem. Of course, these are just some examples of possible applications of the regularity theory for Monge-Amp\`ere.  For instance, as described in the survey paper \cite[Sections 5 and 6]{TWsurvey}, Monge-Amp\`ere equations play a crucial role in affine geometry, for instance in the study of affine spheres and affine maximal surfaces.

\section{Classical theory}

In this section we give a brief overview of some relevant results on the Monge-Amp\`ere equation. Before entering into the concept of weak solutions and their regularity, we first discuss convexity of solutions and the terminology ``degenerate ellipticity'' associated to this equation.

\subsection{On the degenerate ellipticity of the Monge-Amp\`ere equation}
\label{sect:deg ellipt}
{Let $u:\Omega \to \R$ be a smooth solution of \eqref{eq:MAclassical} with $f=f(x)>0$ smooth.
A standard technique to prove regularity of solutions to nonlinear PDEs consists in differentiating the equation solved by $u$ to obtain a linear second-order equation for its
first derivatives. More precisely, let us fix a direction $e \in \mathbb S^{n-1}$ and differentiate \eqref{eq:MAclassical} in the direction $e$. Then,
using the formula
$$
\frac{d}{d\e}\Big|_{\e=0} \det(A+\e B)=\det(A)\,{\rm tr}(A^{-1}B)
\qquad \text{$\forall\,A,B\in\R^{n\times n}$ with $A$ invertible},
$$
we obtain the equation
\begin{equation}
\label{eq:linearMA2}
\det( D^2 u ) \,u^{ij}\p_{ij} u_e=f_e \qquad \text{in $\Omega$}.
\end{equation}
Here $u^{ij}$ denotes the inverse matrix of $u_{ij}:=(D^2u)_{ij}$, lower indices 
denotes partial derivatives  (thus $u_e:=\p_eu)$, and we are summing over repeated indices. Since $\det D^2 u =f>0$, the above equation can be rewritten as
\begin{equation}
\label{eq:linearMA}
a_{ij} \p_{ij} u_e=\frac{f_e}{f} \qquad \text{in $\Omega$},\qquad \text{where } a_{ij}:=u^{ij}.
\end{equation}
Thus, to obtain some regularity estimates on $u_e$, we would like the matrix $a_{ij}$
to be positive definite in order  to apply elliptic regularity theory for linear equations.
But for the matrix $a_{ij}=u^{ij}$ to be positive definite we need $D^2u$ to be positive definite,
which is exactly the convexity assumption on $u$.\footnote{Of course the theory would be similar if one assumes $u$ to be concave. The real difference arises if the Hessian of $u$ is indefinite, since \eqref{eq:linearMA} becomes hyperbolic (and the equation is then called ``hyperbolic Monge-Amp\`ere''). This is still a very interesting problem, but the theory for such equation is completely different from the one of the classical Monge-Amp\`ere equation and it would go beyond the scope of this note.}

We also observe that, without any a priori bound on $D^2u$, the matrix $a_{ij}$ may have arbitrarily small eigenvalues
and this is why one says that \eqref{eq:MAclassical} is ``degenerate elliptic''.
\smallskip

Notice that if one can show that \begin{equation}
\label{eq:C11 both}
c_0{\rm Id} \leq D^2u \leq C_0{\rm Id} \qquad \text{inside $\Omega$}
\end{equation}
for some positive constants $c_0,C_0>0$,
then $C_0^{-1}{\rm Id} \leq (a_{ij})_{1\leq i,j\leq n} \leq c_0^{-1}{\rm Id}$ and the linearized equation \eqref{eq:linearMA}
becomes uniformly elliptic. For this reason, proving \eqref{eq:C11 both}
is one of the key steps for the regularity of solutions to \eqref{eq:MAclassical}.

In this regard we observe that, under the assumption $f(x)\geq \lambda>0$, the product of the eigenvalues of $D^2u$ (which are positive) is bounded from below. Thus, if one can prove that
$|D^2u|\leq C$, one easily concludes that \eqref{eq:C11 both} holds (see \cite[Remark 1.1]{figalliBook} for more details).

In conclusion,  the key step towards the smoothness of solutions consists in proving that $D^2u$ is bounded.}

\subsection{Alexandrov solutions}
\label{sect:Alex sol}

 In his study of the Minkowski problem,
Alexandrov introduced a notion of weak solution to the Monge-Amp\`ere equation that allowed him to give a meaning to the Gaussian curvature of non-smooth convex sets. We now introduce this fundamental concept.

Given an open convex domain $\Omega$,
 the subdifferential of a convex function \(u:\Omega \to \R\) is given by
 \[
\partial u (x):=\{p\in \R^{n}\,:\, u(y)\ge u(x)+p\cdot (y-x)\quad \forall\, y \in \Omega\}.
\]
One then defines the \emph{Monge-Amp\`ere measure} of \(u\) as follows:
$$
\mu_{u}(E):=|\partial u (E)| \qquad \text{for every Borel set $E\subset \Omega$,}
$$
where 
$$
\partial u (E):=\bigcup_{x\in E}\partial u(x)
$$
and $|\cdot|$ denotes the Lebesgue measure.
It is possible to  show that \(\mu_{u}\) is a Borel measure (see \cite[Theorem 2.3]{figalliBook}). Note that,
in the case \(u\in C^{2}(\Omega)\), the change of variable formula gives
$$
|\partial u (E)|=|\nabla u(E)|=\int_E \det D^{2} u (x)\,dx\qquad \text{for every Borel set $E\subset \Omega$,}
$$
therefore
\[
\mu_{u}=\det D^{2} u(x)\,dx
\]
(see \cite[Example 2.2]{figalliBook}).

This discussion motivates the following definition:

\begin{defi}[Alexandrov solutions]\label{ch1:def:alesol}Given an open convex set \(\Omega\) and a function $f:\Omega\times \R\times \R^n\to \R^+$, a convex
function \(u:\Omega \to \R\) is called an \emph{Alexandrov solution} to the Monge-Amp\`ere equation
\[
\det D^{2} u =f(x,u,\nabla u)\qquad \text{in $\Omega$}
\]
if \(\mu_{u}=f(x,u,\nabla u)\,dx\) as Borel measures, namely
$$
\mu_u(A)=\int_Af(x,u,\nabla u)\,dx\qquad \forall\,A\subset \Omega \text{ Borel}.
$$
\end{defi}

Note that because convex functions are locally Lipschitz, they are differentiable a.e. Thus $f(x,u,\nabla u)$ is defined a.e. and the above definition makes sense.

To simplify the presentation, we shall discuss only the case $f=f(x)$, although all the arguments can be extended to the case $f=f(x,u,\nabla u)$ under the assumption that $\partial_u f\geq 0$ (this is needed to ensure that the maximum principle holds, see \cite[Chapter 17]{GT}).

Actually, even if one is interested in solving the Monge-Amp\`ere equation with a smooth right hand side, in order to prove existence of solutions it will be useful to consider also Borel measures as right hand sides.
So, given a nonnegative Borel measure $\nu$ inside $\Omega,$ we shall say that $u$ is an Alexandrov solution of $\det D^2u=\nu$ if $\mu_u=\nu$. 
\smallskip

A fundamental property of the Monge-Amp\`ere measure is that it is stable under uniform convergence
(see \cite[Proposition 2.6]{figalliBook}):

\begin{prop}\label{uni}
Let $u_k:\Omega \to \R$ be a sequence of convex functions converging locally uniformly to $u$.
Then the associated Monge-Amp\`ere measures $\mu_{u_k}$ weakly$^\ast$
converge to $\mu_u$, i.e.
$$
\int_\Omega \varphi\,d\mu_k \to \int_\Omega \varphi\,d\mu\qquad \forall\,\varphi\in C_c(\Omega).
$$
\end{prop}

Another crucial property of this definition is the validity of a comparison principle (see \cite[Theorem 2.10]{figalliBook}):

\begin{prop}\label{thm:comparison}
Let $\mathcal U\subset \Omega$ be an open bounded set, and let $u,v:\Omega\to \R$ be two convex functions satisfying
$$
\begin{cases}
\mu_u \leq \mu_v &\text{in }\mathcal U\\
u\geq v &\text{on }\partial \mathcal U.
\end{cases}
$$
Then
$$
u\ge v
\qquad \text{in }\mathcal U.
$$
\end{prop}

A direct consequence of this result is the uniqueness and stability of solutions (see \cite[Corollaries 2.11 and 2.12]{figalliBook}):

\begin{coro}\label{cor:uniqueness}
Let $\Omega$ be an open bounded set, and $\nu_k:\Omega\to \R$ a family of nonnegative Borel measures satisfying $\sup_k\nu_k(\Omega)<\infty$.
Then, for any $k$, there exists at most one convex function $u_k:\Omega\to \R$ solving the Dirichlet problem
$$
\begin{cases}
\mu_{u_k} =\nu_k &\text{in }\O\\
u_k=0 &\text{on }\partial \O.
\end{cases}
$$
In addition, if $\nu_k\rightharpoonup^*\nu_\infty$ and the solutions $u_k$ exist, then $u_k\to u_\infty$ locally uniformly, where $u_\infty$ is the unique solution of
$$
\begin{cases}
\mu_{u_\infty} =\nu_\infty &\text{in }\O\\
u_\infty=0 &\text{on }\partial \O.
\end{cases}
$$
\end{coro}

Finally, exploiting these results, one can prove existence of solutions (see \cite[Theorem 2.13]{figalliBook}):

\begin{theo}
\label{thm:existence zero}
Let $\Omega$ be an open bounded convex set, and let $\nu$ be a nonnegative Borel measure with $\nu(\Omega)<\infty$.
Then there exists a unique convex function $u:\O\to \R$ solving the Dirichlet problem
\begin{equation}
\label{eq:Dir 0 bis}
\begin{cases}
\mu_u =\nu &\text{in }\O\\
u=0 &\text{on }\partial \O.
\end{cases}
\end{equation}
\end{theo}
\begin{proof}[Sketch of the proof]
Since uniqueness follows from Corollary~\ref{cor:uniqueness}, one only needs to show existence.

By the stability in Corollary~\ref{cor:uniqueness},
since any finite measure can be approximated in the weak$^\ast$ topology by a finite sum of Dirac deltas, 
it suffices to solve the Dirichlet problem \eqref{eq:Dir 0 bis} when $\nu=\sum_{i=1}^N \a_i\delta_{x_i}$ with $x_i \in \O$ and $\a_i>0$.

To prove existence of a solution, one uses the so-called ``Perron Method'': one defines 
$$
\mathcal S[\nu]:=\{v:\O\to\R \text{ convex}\,:\,v|_{\p\O}=0,\,\mu_v \geq \nu \text{ in }\O\},
$$
and shows that this set is nonempty and that it is closed under maximum (namely, $v_1,v_2 \in \mathcal S[\nu]\Rightarrow\max\{v_1,v_2\}\in \mathcal S[\nu]$).
Thanks to these properties, one obtains that $u:=\sup_{v\in \mathcal S[\nu]}v$ is still an element of $\mathcal S[\nu]$, and then one exploits the maximality of $u$ to deduce that $\mu_u=\nu$.

We refer to \cite[Proof of Theorem 2.13]{figalliBook} for more details.
\end{proof}
Actually, if $\Omega$ is strictly convex, a similar argument combined with the existence of suitable barriers allows one to prove the existence of solutions for any continuous boundary datum (see for instance \cite[Theorem 2.14]{figalliBook}):
\begin{theo}
\label{thm:existence g}
Let $\Omega$ be an open bounded strictly convex set, let $\nu$ be a nonnegative Borel measure with $\nu(\Omega)<\infty$, and let $g:\partial\Omega\to \R$ be a continuous function.
Then there exists a unique convex function $u:\O\to \R$ solving the Dirichlet problem
\begin{equation}
\label{eq:Dir g bis}
\begin{cases}
\mu_u =\nu &\text{in }\O\\
u=g &\text{on }\partial \O.
\end{cases}
\end{equation}
\end{theo}

\subsection{Existence of smooth solutions and global regularity}
\label{sect:continuity}

As shown in the previous section, uniqueness of solutions to the Dirichlet problem
holds even at the level of weak solutions.
So, the main issue is existence.

Existence of smooth solutions to the Monge-Amp\`ere equation dates back to the work of Pogorelov \cite{pog smooth}. This is obtained
through the well-celebrated \emph{method of continuity} 
that we now briefly describe (we refer to \cite[Chapter 17]{GT} and \cite[Section 3.1]{figalliBook} for a more detailed exposition).
\smallskip

Assume that $\Omega$ is a smooth uniformly convex domain,\footnote{We say that a domain is uniformly convex if  there exists a radius \(R\) such that
\[\Omega\subset B_{R}(x_0+R\nu_{x_0})
\qquad \text{for every \(x_0\in \partial \Omega\)},
\] where \(\nu_{x_0}\) is the interior normal to \(\Omega\) at \(x_0\). Note that, for a smooth  domain, this is equivalent to asking that the second fundamental form of \(\partial \Omega\) is uniformly positive definite.} and consider $\bar u:\Omega\to \R$ a smooth uniformly convex function
that vanishes on $\partial\Omega$. Then, if we set $\bar f:=\det D^2\bar u$, we have that $\bar f>0$ in $\overline\Omega$ and $\bar u$ solves
\[
\begin{cases}
\det D^2 \bar u=\bar f\quad&\text{in \(\Omega\)}\\
\bar u =0 &\text{on \(\p \Omega\)}.
 \end{cases}
 \]
Now, assume we want to solve
 \begin{equation}\label{maloc}
\begin{cases}
\det D^2  u= f\quad&\text{in \(\Omega\)}\\
 u =0 &\text{on \(\p \Omega\)}
 \end{cases}
 \end{equation}
 for some given $f:\Omega\to \R$ with $f>0$.
Define \(\{f_t:=(1-t)\bar f+tf\}_{t \in [0,1]}\), \(t\in [0,1]\), and  consider the \(1\)-parameter family of problems
  \begin{equation}\label{MAt}
\begin{cases}
\det D^2  u_t= f_t\quad&\text{in \(\Omega\)}\\
 u_t =0 &\text{on \(\p \Omega\)}.
 \end{cases}
 \end{equation}
The method of continuity consists in showing that  the set of \(t \in [0,1]\) such that \eqref{MAt} is smoothly solvable is both open and closed. Since the problem is solvable for $t=0$ (because $\bar u$ is a solution), this implies the existence  of a smooth solution  to our \eqref{maloc}. 

More precisely, assuming that 
$\Omega$ is a uniformly convex domain of class $C^{2,\alpha}$ for some $\alpha\in(0,1)$, it follows that the function $\bar f={\rm det}D^2\bar u$ belongs to $C^{0,\alpha}(\overline\Omega)$.
Then, assuming that $f \in C^{0,\alpha}(\overline\Omega)$, we can consider the set of functions
 \[
 \mathcal C:=\{v:\overline \Omega \to \R \text{ convex functions of class \(C^{2,\alpha}(\overline \Omega)\), \(v=0\) on \(\p \Omega\)}\},
 \]
 and define the nonlinear map
 \[
 \begin{split}
 \mathcal F\colon \mathcal C\times[0,1]&\longrightarrow C^{0,\alpha}(\overline \Omega)\\
(v,t)\,\,\,\,&\longmapsto  \det D^2 v-f_t.
\end{split}
 \]
The goal is show  that the set
\[
\mathcal T:=\{t\in [0,1]: \text{ there exists a \(u_t \in \mathcal C\) such that \(\mathcal F(u_t,t)=0\)}\},
\] 
is nonempty, and it is both open and closed inside $[0,1]$. We now explain the main steps of the argument.

\smallskip

$\bullet$ Nonemptyness follows from the fact that $\mathcal F(\bar u,0)=0$, thus $0 \in 
\mathcal T$.

\smallskip

$\bullet$ Openness follows from the Implicit Function Theorem in Banach spaces (see \cite[Theorem 17.6]{GT}). Indeed, the Frech\`et differential of \(\mathcal F\) with respect to \(v\) is given by the linearized Monge-Amp\`ere operator (compare with \eqref{eq:linearMA2})
\begin{equation}\label{linearized}
D_u \mathcal F(v,t)[h]=\det( D^2 u )u^{ij} h_{ij},\qquad h=0\, \text{ on \(\p \Omega\),}
\end{equation}
where we set \(h_{ij}:=\partial_{ij} h\),  \(u^{ij}\) is the inverse of \( u_{ij}:=\partial_{ij}u\), and we are summing over repeated indices. Notice that if a function \(v\) is bounded in  \(C^{2,\alpha}\) and  \(\det D^2v\) is bounded from below, then the smallest eigenvalue of \(D^2 v\) is bounded uniformly away from zero and the linearized operator  becomes uniformly elliptic with \(C^{0,\alpha}\) coefficients
(cp. Section \ref{sect:deg ellipt}).
Therefore, classical Schauder's theory gives  the invertibility of  \(D_u \mathcal F(u_t,t)\)  whenever $u_t$ solves $\mathcal F(u_t,t)=0$ (see for instance \cite[Chapter 6]{GT}).

\smallskip

$\bullet$ The proof of closedness is done via global a priori estimates. More precisely, 
the following fundamental a priori bound holds (see \cite[Theorem 3.2]{figalliBook}):\footnote{The assumption \(u\in C^4(\Omega)\) in Theorem \ref{po2} is not essential, as it is needed only to justify the computations in the proof.}
\begin{theo}\label{po2}Let \(\Omega\) be a uniformly convex domain 
of class
\(C^3\), and let \(u\in C^4(\Omega)\) be a solution of \eqref{maloc} with \(f\in C^2(\overline \Omega)\) and  \(0<\lambda \le f\le1/\lambda\). Then there exists a constant \(C\), depending only on \(\Omega\), \(\lambda\), \(\|f\|_{C^2(\overline \Omega)}\), such that
\[
\|D^2 u\|_{C^0(\overline \Omega)}\le C.
\] 
\end{theo}
 As already noticed in Section \ref{sect:deg ellipt}, once a uniform bound on \(D^2u\) inside \(\overline \Omega\) holds, the Monge-Amp\`ere equation becomes uniformly elliptic and
classical elliptic regularity theory  yields \(C^{2,\alpha}\) estimates for solutions of $\mathcal F(u_t,t)=0$, proving the desired closedness of \(\mathcal T\).

\smallskip

Thanks to this argument, one concludes the validity of the following existence result:
\begin{theo}\label{existencesmooth}Let \(\Omega\) be a uniformly convex domain
of class $C^3$. Then, for all \(f\in C^{2}(\overline \Omega)\) with \(0<\lambda \le f\le 1/\lambda\),
there exists a unique \(u \in C^{2,\alpha}(\overline \Omega)\) solution to \eqref{maloc}.
\end{theo}
Recalling that uniqueness holds also at the level of Alexandrov solutions, this proves the $C^{2,\alpha}$ regularity (for any $\alpha<1$) of Alexandrov solutions in $C^3$ uniformly convex domains with $C^2$ right hand side.
It is interesting to remark that the $C^3$ regularity assumption on the boundary is necessary, as shown by Wang in \cite{wang2}.

\subsection{Caffarelli's regularity theory}
\label{sect:Caff}
We now investigate the regularity of Alexandrov solutions under weaker smoothness assumptions on the right hand side.

In the 90's
Caffarelli developed a regularity theory for Alexandrov solutions, showing that strictly convex solutions of \eqref{eq:MAclassical}
are locally $C^{1,\gamma}$ provided $\lambda \leq f\leq 1/\lambda$
for some $\lambda >0$ \cite{Caf1,Caf2,Caf3}. 
We emphasize that, for weak solutions, strict convexity is not implied by the positivity of $f$ (unless $n=2$)
and it is actually necessary for regularity, see Section \ref{sect:strict} below.

The following result is proved in \cite{Caf2}:
\begin{theo}\label{thm:calfa}
Let \(u:\Omega \to \R\) be a {strictly convex} Alexandrov solution of $\mu_u=f\,dx$ with \(0<\lambda \le f\le 1/\lambda\). Then \(u\in C_{\rm loc}^{1,\gamma}(\Omega)\) for some  \(\gamma=\gamma(n,\lambda)>0\).
\end{theo}

To explain the idea behind the proof of the above theorem, let us point out the following simple properties of solutions to the Monge-Amp\`ere equation (this is another manifestation of its degenerate ellipticity):
If \(A:\R^n\to \R^n\) is an affine transformation with \(\det A=1\),\footnote{Given an affine transformation $Ax:=Mx+v$, by abuse of notation we write $\det A$ in place of $\det M$.} and \(u\) is a solution of the Monge-Amp\`ere equation
with right hand side $f$, then \(u\circ A\) is a solution to the Monge-Amp\`ere equation with right hand side $f \circ A$. This affine invariance creates serious obstructions to obtain a local regularity theory.
Indeed, for instance, the  functions  
 \[
 u_\e(x_1,x_2)=\frac{ \e x_1^{2}}{2}+\frac{x_2^{2}}{2\e}-1
 \]
 are solutions to \(\det D^2 u_\e=1\) inside the convex set \(\{u_\e < 0\}\). Thus, unless the level set $\{u_\e=0\}$ is sufficiently ``round'', there is no hope to obtain a priori estimates on \(u\). The intuition of Caffarelli was to use the so-called John's Lemma \cite{John}:
 
\begin{lemm}
\label{lem:john} Let \(\mathcal K \subset \R^n\) be a bounded convex set with non-empty interior. Then there exists an ellipsoid \(E\) satisfying
\begin{equation}\label{ap:j}
E\subset \mathcal K\subset nE,
\end{equation}
where $nE$ denotes the dilation of $E$ by a factor $n$ with respect to its center. 
\end{lemm}
In particular, if we define a convex set \(\mathcal K\)  to be \emph{normalized} if 
\[
B_1\subset \mathcal K \subset nB_{1},
\]
then Lemma \ref{lem:john} says that, for every bounded open convex set \(\mathcal K\), there is an affine transformation \(A:\R^n\to \R^n\) such that \(A(\mathcal K)\) is normalized. 

Note that, if $u$ is strictly convex,
given a point $x\in \Omega$ and $p \in \partial u(x)$ one can choose $t>0$ small enough so that the convex set
\begin{equation}
\label{sec}
S(x,p,t):=\{z\in \Omega\,:\,u(z)-u(x)-p\cdot (z-x)<t\}
\end{equation}
is contained inside $\Omega$.
Then, if we replace $u(z)$ with $u_x(z):=u(z)-u(x)-p\cdot (z-x)-t$, it follows that
\begin{equation*}
\lambda\,dx \le \mu_{u_x}\le \frac{1}\lambda\,dx\quad \text{in \(S(x,p,t)\)},\qquad u_x=0 \quad \text{on \(\partial \bigl(S(x,p,t)\bigr)\)}.
\end{equation*} 
Also, if \(A:\R^n\to \R^n\) normalizes \(S(x,p,t)\), then \(v:=(\det A)^{2/n}\, u_x\circ A^{-1}\) solves 
\begin{equation}
\label{eq:v normalized}
\lambda\,dx \le  \mu_v \le \frac{1}\lambda\,dx\quad \text{in \(A(S(x,p,t))\)},\qquad v=0 \quad \text{on \(\partial \bigl(A(S(x,p,t))\bigr)\)}.
\end{equation}
Thanks to the above discussion, it suffices to prove the result when $u$ is a solution inside a normalized convex set. In other words, Theorem \ref{thm:calfa} is a direct consequence of the following result:

\begin{theo}\label{thm:calfa2}
Let  \(\Omega\) be a normalized convex set, and \(u\) be a solution of
$$
\mu_u=f\,dx
\quad \text{in \(\Omega\)},\qquad u=0 \quad \text{on \(\partial \Omega\)},
$$
with \(0<\lambda \le f\le 1/\lambda\). Then $u$ is strictly convex inside $\Omega$, and \(u\in C_{\rm loc}^{1,\gamma}(\Omega)\) for some  \(\gamma=\gamma(n,\lambda)>0\).
\end{theo}

In the proof of the above theorem, a key step consists in showing that solutions of \eqref{maloc} inside normalized domains have a \emph{universal} modulus of strict convexity. A fundamental ingredient to prove this fact is the  following important result of Caffarelli \cite{Caf1}
(see also \cite[Theorem 4.10]{figalliBook}):
\begin{prop}\label{ch2:thm:contact}Let \(u\) be a solution of
\[
\lambda\,dx \le  \mu_u \le \frac{1}\lambda\,dx
\]
inside a convex set \(\Omega\), $x \in \Omega,$ and $p\in \partial u(x)$.
Let \(\ell(z):=u(x)+p\cdot(z-x)\). If  the convex set
\[
W:=\{z\in \Omega\,:\, u(z)=\ell(z)\}
\]
contains more than one point, then it cannot have extremal points inside \(\Omega\).
\end{prop}
This statement says that if a solution coincides with one of its supporting planes on more than one point (that is, it is not strictly convex), then the contact set has to cross the domain.
In particular this is not possible if $u|_{\partial\Omega}=0$ (as otherwise $u\equiv 0$), proving that solutions to \eqref{eq:v normalized} are strictly convex.

\begin{proof}
[Sketch of the proof of Theorem \ref{thm:calfa2}]
As mentioned above, Proposition \ref{ch2:thm:contact} implies that $u$ is strictly convex. Also,
by compactness, one can prove that the modulus of strict convexity of $u$ is universal (see \cite[Section 4.2.2]{figalliBook}).

We then apply this information at all scales. More precisely, given any point $x\in \Omega$, $p \in \partial u(x)$, and $t>0$ small, we  consider $u_x(z):=u(z)-u(x)-p\cdot (z-x)-t$. Then, if \(A:\R^n\to \R^n\) normalizes \(S(x,p,t)\), the function \(v:=(\det A)^{2/n}\, u_x\circ A^{-1}\) enjoys the same strict convexity properties as $u$.
Using this fact at all points $x$
and for all small values of $t$, a careful iteration argument proves the validity of Theorem \ref{thm:calfa2} (see the proof of \cite[Theorem 4.20]{figalliBook} for more details).
\end{proof}

\medskip

Note that Theorem \ref{existencesmooth} is unsatisfactory from a PDE viewpoint: indeed, it requires the $C^2$ regularity of $f$ to prove the $C^{2,\alpha}$ regularity of the solution, while the usual elliptic regularity theory would suggest that $f \in C^{0,\alpha}$ should be enough. This is indeed true, as proved 
by Caffarelli in
\cite{Cafw2p} (again, it suffices to consider normalized convex sets):

\begin{theo}
\label{thm:C2a}Let  \(\Omega\) be a normalized convex set, and \(u\) be a solution of
$$
\mu_u=f\,dx
\quad \text{in \(\Omega\)},\qquad u=0 \quad \text{on \(\partial \Omega\)},
$$
with \(0<\lambda\le f\le 1/\lambda\) and  \(f\in C_{\rm loc}^{0,\alpha}(\Omega)\).
Then $u \in C^{2,\alpha}_{\rm loc}(\Omega)$.
\end{theo}

The proof of the above theorem is based on the property that, under the assumption that \(f\) is almost constant (say very close to \(1\)), \(u\) is very close to a solution of $\mu_v=dx$. Since this latter function is smooth (by Theorem \ref{existencesmooth}),
an iteration argument permits to show that the \(C^{2,\alpha}\) norm of \(u\) remains bounded (see also \cite[Theorem 4.42]{figalliBook}).

With this line of reasoning one  can also prove the following theorem \cite{Cafw2p}:

\begin{theo}\label{thm:W2p}
Let  \(\Omega\) be a normalized convex set, and \(u\) be a solution of
$$
\mu_u=f\,dx
\quad \text{in \(\Omega\)},\qquad u=0 \quad \text{on \(\partial \Omega\)}.
$$
Then, for every \(p>1\) there exists a positive constant \(\delta(p)\) such that if \(\|f-1\|_\infty \le \delta(p)\) then $u \in W^{2,p}_{\rm loc}(\Omega)$.
\end{theo}

Since any continuous function is arbitrarily close to a constant at small scales, one obtains the following:

\begin{coro}
Let  \(\Omega\) be a normalized convex set, and \(u\) be a solution of
$$
\mu_u=f\,dx
\quad \text{in \(\Omega\)},\qquad u=0 \quad \text{on \(\partial \Omega\)},
$$
with $f>0$ continuous. Then
$u \in W^{2,p}_{\rm loc}(\Omega)$ for any $p<\infty$.
\end{coro}

\begin{rema}
\label{rmk:W2p}
As shown in \cite{F-proc}, exploiting the ideas introduced in \cite{DF1,DFS} one can find an explicit estimate for $\delta(p)$ in terms of $p$ in Theorem \ref{thm:W2p}, namely
$\delta(p)\simeq e^{-C\,p}$ for some dimensional constant $C>0$.
\end{rema}

\subsection{Some applications}
\label{sect:applic classic}

In this section we briefly describe two applications of the regularity theory developed in the previous sections.

\subsubsection{The Minkowski problem}

A classical problem in convex geometry is to prescribe some geometric quantity (the surface area, the Gaussian curvature, etc.) and find necessary and sufficient conditions ensuring that such a quantity comes from a convex domain.
In this section we briefly discuss the ``prescribed Gaussian curvature'' problem.

Let $K\subset \R^n$ be an open bounded convex domain containing the origin, and parameterize $\p K$ in polar coordinates as follows:
$$
\p K=\bigl\{\rho(x)\,x\,:\,x \in \mathbb S^{n-1},\, \rho:\mathbb S^{n-1}\to \R^+\bigr\}.
$$
Then, to any point $z \in \p K$ we associated the {normal mapping}
$$
N_K(z):=\Bigl\{y \in \mathbb S^{n-1}\,:\, K\subset \{y\,:\,\<y,w-z\>\leq 0\}\Bigr\}.
$$
Geometrically, the normal mapping finds the normals of all supporting hyperplanes at $z$,
and we can think of $N_K$ as an analogue of the subdifferential map.

Finally, we consider the (multivalued) {Gauss map} $G_K:\mathbb S^{n-1}\to \mathbb S^{n-1}$ defined by
$$
G_K(x):=N_K\bigl(\rho(x)\,x\bigr),
$$
and define the {Gaussian curvature measure}
$$
\mu_K(E):=\HH^{n-1}\bigl(G_K(E)\bigr)\qquad \forall\,E\subset \mathbb S^{n-1}\text{ Borel},
$$
where $\HH^{n-1}$ denotes the $(n-1)$-dimensional Hausdorff measure on $\mathbb S^{n-1}$. 

As for the Monge-Amp\`ere measure, one can show that $\mu_K$ is a Borel measure. 
 One then asks the following question: {\it Given a Borel measure $\nu$ on $\mathbb S^{n-1}$, can we find an open bounded convex set $K$ containing the origin and such that $\mu_K=\nu$?}

In \cite{AL0,AL1}, Alexandrov found necessary and sufficient conditions to ensure the existence of a solution to this problem. 
As for the existence of Alexandrov solutions to the Monge-Amp\`ere equation,
the existence of $K$ is first proved when $\nu$ is a finite sum of Dirac deltas, and  
then one obtains the general case by approximation.
The original existence proof of Alexandrov when $\nu$ is discrete was based on a topological argument
relying on the Invariance of Domain Theorem \cite{AL0}
(see also \cite{AL3}).

Thanks to the regularity theory developed by Caffarelli, one obtains the following regularity result:
\begin{theo}
Let $K\subset \R^n$ be an open bounded convex domain containing the
origin, and assume that $\mu_K = f\,d\mathcal H^{n-1}$ for some $f : \mathbb S^{n-1} \to \mathbb R$, with $0<\lambda \leq f \leq 1/\l$. Then $\partial K$ is strictly convex and  of class
$C^{1,\gamma}$. If in addition $f \in C^{0,\alpha}$ for some $\alpha\in (0,1)$, then $\partial K \in C^{2,\alpha}$.
\end{theo}	

\begin{proof}[Sketch of the proof]
Since $K$ is convex, one can locally parameterize the boundary as the graph of a convex function $u:\Omega \subset \R^{n-1}\to \R$ .
It is a classical fact that the Gaussian curvature of the graph
of a $C^2$ function $v:\Omega \subset \R^{n-1}\to \R$ 
is given by
$$
\frac{\det D^2 v}{(1+|\n v|^2)^{\frac{n+1}2}}.
$$
Then, by the assumption $\mu_K = f\,d\mathcal H^{n-1}$, an approximation argument based on Proposition~\ref{uni} yields the validity of the equation
$$
\mu_u=f(x)\bigl(1+|\n u|^2\bigr)^{\frac{n+1}2}\,dx,
$$
where $\n u$ exists at almost every point since $u$ is locally Lipschitz (being convex).
In particular, $\mu_u$ is locally bounded. Applying Proposition \ref{ch2:thm:contact} and Theorem \ref{thm:calfa}, one deduces that $\partial K$ is strictly convex and of class $C^{1,\gamma}$. Finally, the $C^{2,\alpha}$ regularity when $f \in C^{0,\alpha}$ follows from Theorem \ref{thm:C2a}.
\end{proof}

\begin{rema}
The regularity theory for Monge-Amp\`ere plays a crucial role  in many other variants of the Minkowski problem. For instance, it appears in the proof of existence and uniqueness of convex domains with prescribed harmonic measure \cite{jerison}.
\end{rema}

\subsubsection{The optimal transport problem}

Let $\mu$ and $\nu$ denote two probability measures on $\R^n$.
The {optimal transport problem} (with quadratic cost) consists in finding the ``optimal'' way of transporting
$\mu$ onto $\nu$ given that the transportation cost to move a point from $x$ to $y$ is $|x-y|^2.$
Hence, one is naturally led to minimize
$$
\int_{\R^n} |S(x)-x|^2\,d\mu(x)
$$
among all maps $S$ that ``transport $\mu$ onto $\nu$''. Mathematically, this corresponds to saying that $S_\#\mu=\nu$, that is, for any bounded Borel function $\varphi:\R^n\to \R$,
$$
\int_{\R^n} \varphi(y)\,d\nu(y)=\int_{\R^n}\varphi\bigl(S(x)\bigr)\, d\mu(x).
$$
By a classical theorem of Brenier \cite{bre} (see also \cite{CAM,RR,McCann1995}),
existence and uniqueness of optimal maps hold provided that $\mu$ is absolutely continuous. Moreover, such a map is given by the gradient of a convex function.
This is summarized in the next theorem:
\begin{theo}
\label{thm:Brenier}
Let $\mu,\nu$ be probability measures on $\R^n$ with $\mu=f\,dx$ and $\nu=g\,dy$.
Then:
\begin{enumerate}
\item[-] There exists a $\mu$-a.e. unique optimal transport map $T$.
\item[-] There exists a lower semicontinuous convex function $u:\R^n\to \R\cup\{+\infty\}$ such that
$T=\nabla u$ $\mu$-a.e.
and
$$
{\rm det}(D^2 u)=\frac{f}{g\circ\nabla u}\qquad \text{$\mu$-a.e.}
$$
\end{enumerate}
\end{theo}
The above theorem shows that optimal transport maps solve a Monge-Amp\`ere equation in a weak sense, that is usually referred to as ``Brenier sense''.

While for Alexandrov solutions one may apply the regularity theory developed in the previous sections, 
 Caffarelli observed in \cite{Caf3} that even for smooth densities
one cannot expect any general regularity result for Brenier solutions without
making some geometric assumptions on the support of the target
measure. Indeed, let $n=2$ and suppose that $X=B_1$ is the unit ball centered
at the origin and $Y=\bigl(B_1^+ + \,e_1\bigr) \cup \bigl(B_1^- -
\,e_1\bigr)$ is the union of two half-balls, where
$$
B_1^+:=\bigl(B_1 \cap \{x_1
>0\}\bigr),\qquad B_1^-:=\bigl(B_1 \cap \{x_1 <0\}\bigr),
$$
and $(e_1,e_2)$ denotes the canonical basis of $\R^2$.
Then if $f=\frac{1}{|X|}\mathbf 1_X$ and $g=\frac{1}{|Y|}\mathbf 1_Y$,
the optimal transport map is given by
$$
T(x):=\left\{ \begin{array}{ll} x+\,e_1& \text{if } x_1 >0\\
x-\,e_1 &\text{if }x_1<0,\end{array}\right.
$$
which corresponds to the gradient of the convex function
$u(x)= |x|^2/2+|x_1|$.

Thus, in
order to hope for a regularity result for $u$ we need at least
to assume the connectedness of $Y$.
However, starting from the above
construction and considering a sequence of domains $Y_\e$ where
one adds a small strip of width $\e>0$ to glue together
$\bigl(B_1^+ + \,e_1\bigr) \cup \bigl(B_1^- - \,e_1\bigr)$, one can
also show that for $\e>0$ small enough the optimal map will still
be discontinuous (see \cite{Caf3,F-partial}). Hence, connectedness is not enough to ensure regularity.
As shown by Caffarelli \cite{Caf3,Caf4},  convexity of $Y$ is the right assumption to ensure that a Brenier solution is also an Alexandrov solution, so that the general regularity theory from the previous sections  apply (see also \cite{delanoe,urbas}):

\begin{theo}
\label{thm:regularity transport}
Let $X,Y\subset \R^n$ be two bounded open sets, let
$f,g:\R^n\to \R^+$ be two
probability densities that are zero outside $X,Y$ and are
bounded away from zero and infinity on $X,Y$, respectively.
Denote by  $T=\nabla u:X\to Y$ the optimal transport map provided by Theorem~\ref{thm:Brenier},
and assume that $Y$ is convex.
Then there exists $\gamma>0$ such that $T \in C^{0,\gamma}_{\rm loc}(X)$.
Furthermore, if $f \in C^{k,\alpha}(\overline X)$ and $g\in C^{k,\alpha}(\overline Y)$
for some integer $k\geq 0$ and some $\alpha\in (0,1)$, and if both $X$ and $Y$ are smooth and uniformly convex, then $T:\overline X\to\overline Y$ is a global diffeomorphism of class $C^{k+1,\alpha}$.
\end{theo}

As shown for instance in \cite{FRV}, the convexity of the target is necessary for the continuity of the optimal transport map. Even worse, as recently shown in \cite{yash}, even with constant densities one can construct a discontinuous optimal transport map from a smooth convex domain to a small Lipschitz deformation of itself. 

All these facts motivate the following very natural question: {\it What can one say when the convexity assumption on the target is removed?}
As shown in \cite{FK-partial,DFpartial} (see also \cite{F-partial} for a more precise description of the singular set in two dimensions, and \cite{GO} for a recent variational proof of the result in \cite{FK-partial}),
one can always prove that the optimal transport map is smooth outside
a closed set of measure zero.

\section{Recent developments I: Interior regularity}

In \cite{wang} Wang showed that for any $p>1$ there exists a function $f$
satisfying $0<\lambda(p) \leq f\leq 1/\lambda(p)$ such that $u \not \in W^{2,p}_{\rm loc}$.
This counterexample shows that the results of Caffarelli are more or less optimal.
However, an important question which remained open was whether strictly convex solutions of $\mu_u=f\,dx$
with $0<\lambda \leq f\leq 1/\lambda$
could be at least $W^{2,1}_{\rm loc}$,
or even $W^{2,1+\varepsilon}_{\rm loc}$ for some $\varepsilon=\varepsilon(n,\lambda)>0$. The question of \(W_{\rm loc}^{2,1}\) regularity has been recently solved by De Philippis and Figalli in \cite{DF1}. Following the ideas introduced there, the result has been refined to \(u\in W_{\rm loc}^{2,1+\e}\) for some \(\e>0\) (see  \cite{DFS, Shm}).

\begin{theo}\label{w21eps}Let  \(\Omega\) be a normalized convex set, and \(u\) be a solution of
$$
\mu_u=f\,dx
\quad \text{in \(\Omega\)},\qquad u=0 \quad \text{on \(\partial \Omega\)},
$$
 with \(0<\lambda\le f\le 1/\lambda\).
Then there exists \(\e=\e(n,\lambda)>0\) such that
$u \in W^{2,1+\e}_{\rm loc}(\Omega).$ 
\end{theo}

Again, as in Section \ref{sect:Caff}, the previous result holds for strictly convex solutions of $\mu_u=f\,dx$ with $0<\lambda\le f\le 1/\lambda$.

\begin{proof}[Sketch of the proof]
Given $x\in\Omega$ and $t>0$ small,
we consider the family
\(\{S(x,\nabla u(x),t)\}_{x\in \Omega,\, t>0}\)
as defined in \eqref{sec}.
Thinking of $S_t(x):=S(x,\nabla u(x),t)$
as the ``ball centered at $x$ with radius $t$'',
any subdomain $\Omega'\subset\subset \Omega$ endowed with this family of ``balls''
is a space of homogeneous type  in the sense of Coifman
and Weiss, see \cite{CG, GH, AFT}. In particular Stein's Theorem implies that if 
\[
\mathcal M(D^2 u)(x):=\sup_{t>0}\mean{S(x,\nabla u(x),t)} |D^2 u| \in L^1_{\rm loc}(\Omega),
\]
then \(|D^2 u|\in L\log L_{\rm loc}\), that is $\int_{\Omega'} |D^2 u| \log(2+|D^2u|) \leq C(\Omega')$ for any $\Omega'\subset\subset \Omega$. 
The key estimate in \cite{DF1} consists in showing that
\[
\|\mathcal M(D^2 u)\|_{L^1_{\rm loc}(\Omega)}\le C \|D^2 u\|_{L^1_{\rm loc}(\Omega)},
\]
for some constant \(C=C(n,\lambda)\).

Once this estimate is proved, it follows by the convexity of $u$ that $L^1_{\rm loc}$ norm of $D^2u$
is locally bounded (see \cite[Equation (4.74)]{figalliBook}), thus\footnote{Here the reader may be confused by the sentence ``Since $u$ is convex, the $L^1_{\rm loc}$ norm of $D^2u$
is locally bounded''. Indeed, this seems to say that the $W^{2,1}_{\rm loc}$ regularity of $u$ is trivial since the integral of $|D^2u|$ is locally finite. This is not the case because, for a convex function, $D^2u$ may be a measure and so $\int_E|D^2u|$ denotes the integral over a set $E$ of the measure $|D^2u|$. So, to prove that $u \in W^{2,1}_{\rm loc}$ is not enough to show that $\int |D^2u|$ is locally finite but one needs to show that $|D^2u|$ is absolutely continuous with respect to the Lebesgue measure.} 
\begin{equation}
\label{eq:LlogL}
|D^2 u| \log(2+|D^2u|)\in L^1_{\rm loc}(\Omega).
\end{equation}
By this a priori estimate and an approximation argument with smooth solutions, as shown in \cite{DF1} one easily deduce that $D^2u$ is an $L^1$ function, and therefore $u \in W^{2,1}_{\rm loc}$.

We now explain how this argument actually implies that $u \in W^{2,1+\e}_{\rm loc}$.
In view of \eqref{eq:LlogL}, the measure where $|D^2 u|$ is large decay in a quantitative way:
 $$|\{ |D^2 u | \ge M\}| \le \frac{1}{M\log M}\int_{\{ |D^2 u | \ge M\}|}|D^2 u| \log(2+|D^2u|) \leq \frac{C}{M\log M},$$
 for any $M$ large.
 In particular, choosing first $M$ sufficiently large and then taking $\e>0$ small enough,
we deduce (a localized version of) the bound
 $$
 |\{ |D^2 u | \ge M\}| \le \frac{1}{M^{1+2\e}} |\{ |D^2 u | \ge 1\}|
 $$
Applying this estimate at all scales (cp. the sketch of the proof of Theorem \ref{thm:calfa2}) together with a covering lemma yields
$$
 |\{ |D^2 u | \ge M^k\}| \le \frac{1}{M^{(1+2\e)k}} |\{ |D^2 u | \ge 1\}|\qquad \forall\,k\geq 1,
$$
and the local $L^{1+ \e}$ integrability for $|D^2 u|$ follows (see for instance \cite[Section 4.8.4]{figalliBook} for more details).
\end{proof}

\subsubsection{An application: the semigeostrophic equations}
The semigeostrophic equations are a simple model used
in meteorology to describe large scale atmospheric flows,
and can be derived from the 3-d 
Euler equations, with Boussinesq and hydrostatic approximations, subject to a strong
Coriolis force \cite{cu}.
Since for large scale atmospheric flows the Coriolis force dominates
the advection term, the flow is mostly bi-dimensional. For this reason,
the study of the semigeostrophic equations in 2-d or 3-d is pretty similar,
and in order to simplify our presentation we focus here on the 2-dimentional periodic case.

The semigeostrophic system can be written as
\begin{equation}\label{eqn:SGsystem2}
\begin{cases}
\partial_t \nabla p_t + (\uu_t\cdot \nabla) \nabla p_t  +\nabla^\perp p_t +\uu_t = 0\\
\nabla \cdot \uu_t = 0 \\
p_0= \bar p
\end{cases}
\end{equation}
where $\uu_t:\R^2 \to \R^2$ and $p_t:\R^2 \to \R$ are periodic functions
corresponding respectively to the velocity and the pressure, and \(\nabla^\perp p_t \) is the \(\pi/2\) counterclockwise rotation of \(\nabla p\).

As shown in \cite{cu},
energetic considerations show that it
is natural to assume that $p_t$ is ($-1$)-convex, i.e., the function
$P_t(x):=p_t(x)+|x|^2/2$ is convex on $\R^2$.
Let  $P^*_t$ be the convex conjugate of $P_t$, namely
\[
P_t^*(y):=\sup_{x\in\R^2} \bigl\{y \cdot x-P_t(x)\bigr\}.
\]
Then, assuming that $0<\lambda\leq \det(D^2P_0^*)\leq 1/\lambda$,  one can prove that
$$
0<\lambda \leq \det(D^2P_t^*)\leq 1/\lambda \qquad \forall\,t >0
$$
in the Alexandrov sense
(see \cite{ACDF1} for more details).
Thanks to Theorem \ref{w21eps} this implies that
$P_t^* \in W^{2,1+\e}_{\rm loc}$, which is one of the key ingredients to prove
the global existence of distributional solutions to \eqref{eqn:SGsystem2}
on the 2-dimentional torus \cite{ACDF1} and in three dimensional domains \cite{ACDF2}.

\begin{rema}
From a physical point of view, the lower bound on $\det(D^2P_0^*)$ is not natural, and it would be very useful if the $W^{2,1}$ regularity of solutions to 
$\mu_u \leq \frac{1}{\lambda}\,dx $ was true, at least in two dimensions. Unfortunately this is false, as shown by Mooney in \cite{Mooney-counterex}.

On a different direction, one would like to prove global existence of smooth solutions of \eqref{eqn:SGsystem2} when the initial datum is smooth. Motivated by the analogous result for the 2d incompressible Euler equation, a possible strategy to prove this result would be to show that strictly convex solutions of $\mu_u=f\,dx$ with $f \in C^{0,\alpha}$ such that $0<\l\leq f \leq 1/\l$ satisfy
$\|D^2u\|_{C^{0,\alpha}}\leq C(n,\lambda,\alpha)\|f\|_{C^{0,\alpha}}$ (namely, the control is linear with respect to the norm of $f$). As shown in \cite{FJM-C2a} this is false, and the global existence of smooth solutions is still an open problem.
\end{rema}

\subsection{On the strict convexity of weak solutions}
\label{sect:strict}

{
As already mentioned,
strict convexity is not just a technical assumption but it is necessary to obtain regularity. Indeed, as discovered by Pogorelov, there exist Alexandrov solutions to the Monge-Amp\`ere equation with smooth positive right-hand side which are not $C^2$.
For instance, the function
\begin{equation}
\label{eq:example pogorelov}
u(x_1,x'):=|x'|^{2-2/n} (1+x_1^2),\qquad(x_1,x')\in \R\times \R^{n-1},\qquad n\ge 3,
\end{equation}
is $C^{1,1-2/n}$ and solves $\det D^2u=c_n(1+x_1^2)^{n-2}(1-x_1^2)>0$ inside $B_{1/2}$.
Furthermore, the bound \(0<\lambda \le \det D^2 u\le 1/\lambda\) is not even enough for $C^1$
regularity:
the function
\[
u(x_1,x'):=|x'|+|x'|^{n/2}(1+x_1^2),\qquad(x_1,x')\in \R^n,\qquad n\ge 3,
\]
is merely Lipschitz and solves \(0<\lambda \le \det D^2 u\le 1/\lambda\)
in a small convex neighborhood of the origin.\footnote{Actually, for $n\geq 3$, one can even construct a Lipschitz Alexandrov solution of $\det D^2u=1$ in a small ball $B_\rho(0)$. To see this, let $\eta>0$ and set $v_\eta(x):=\eta\left(|x'|+|x'|^{n/2}(1+x_1^2)\right).$ Then, if $\eta>0$ is large enough, it follows that $\det D^2v_\eta\geq 1$ inside $B_\rho(0)$ for some $\rho>0$ small.

Let $w_\eta:B_\rho(0)\to \R$ be the convex envelope of $v_\eta|_{\partial B_\rho(0)}$. It is a classical fact that $\det D^2w_\eta=0$ in the Alexandrov sense (see for instance \cite{OS}). Also, since $v_\eta\geq 0$ it follows that $w_\eta\geq 0$.
Finally, since $v_\eta(x_1,0)=0$ for $x_1=\pm \rho$, $w_\eta(x_1,0)=0$ for $|x_1|\leq \rho$.

Now, let $u$ be the Alexandrov solution of
$$
\begin{cases}
\det D^2u =1 &\text{in }B_\rho\\
u=v_\eta &\text{on }\partial B_\rho
\end{cases}
$$
provided by Theorem \ref{thm:existence g}.
Then it follows by Proposition \ref{thm:comparison} that $v_\eta\leq u\leq w_\eta$ inside $B_\rho$.
This implies in particular that $u(x_1,0)=0$ for $|x_1|\leq \rho$, that combined with
$$
u(x_1,x')\geq v_\eta(x_1,x')\geq \eta |x'|
$$
shows that $u$ is merely Lipschitz continuous.
}

Alexandrov showed in \cite{AL2} that, in contrast with the above counterexamples, every two dimensional solution of  \(\mu_u \ge \lambda\,dx>0\) is strictly convex. In \cite{Caf-CPDE}, Caffarelli generalized these examples to solutions that degenerate along
subspaces, and he proved that solutions can degenerate only on subspaces of dimension less than $n/2$.}

Since one cannot hope for $C^1$
regularity  of non-strictly convex solutions, it is natural to ask whether one can obtain some integrability estimates for the second derivatives.
In the previous section  we showed that strictly convex solutions of $0<\l\,dx \leq \mu_u \leq \frac1\l\,dx$
are $W^{2,1+\e}_{\rm loc}$ for some $\e=\e(n,\l)>0$.
If one denotes by $\Sigma$ the ``singular set'' of points where $u$ is not strictly convex, that is
$$
\Sigma:=\{x \in \O\,:\, \exists \,z \in \O\setminus\{x\} \text{ and }p\in\p u(x)\quad \text{s.t.} \quad u(z)=u(x)+\<p,z-x\> \},
$$
then one may wonder whether the second derivatives of $u$ can concentrate on $\Sigma$.
This fact has been recently ruled out by Mooney \cite{Mooney}
who showed that the $(n-1)$-dimensional Hausdorff measure of  $\Sigma$ vanishes.
From this, he deduced the $W^{2,1}$ regularity
of solutions without any strict convexity assumptions.
Actually, in a subsequent paper \cite{MooneyW21}, he was able to strengthen this result by showing a small logarithmic integrality improvement
and proving that such a result is optimal.

\begin{theo}
Let $\O\subset \R^n$ be an open set, and $u:\O\to \R$ be a convex function satisfying $\mu_u=f\,dx$ for some $0<\l\leq f \leq 1/\l$.
Then $\HH^{n-1}(\Sigma)=0$ and $u \in W^{2,1}_{\rm loc}(\O)$.
In addition, there exists $\eta=\eta(n)>0$ such that
$$
\int_{\O'}|D^2 u|\,\log\bigl(2+|D^2u|\bigr)^\eta\,dx <\infty\qquad \forall\,\O'\subset \subset \Omega.
$$
On the other hand, if  $M>0$ is sufficiently large,  one can construct a solution $u$ with $f\equiv 1$ such that 
$$
\int_{\O'}|D^2 u|\,\log\bigl(2+|D^2u|\bigr)^M\,dx =+\infty\qquad \text{for some }\Omega'\subset\subset \O.
$$
\end{theo}

\section{Recent developments II: Boundary regularity}
\label{sect:bdry reg}

The interior regularity theory for Alexandrov solutions relies on several geometric properties of sections $S(x,p,y)$ of $u$ that are strictly contained inside $\O$ (see \eqref{sec}).
In particular, an important property states that any ``section'' $S_t(x):=S(x,\n u(x),t)$ 
contained inside $\O$ is comparable to an ellipsoid of volume $t^{n/2}$ (see for instance \cite[Lemma 4.6]{figalliBook}).

In order to develop a boundary regularity theory, it is crucial to understand the geometry of sections $S_t(x)$
when $x\in \p\O$. This has been done by Savin in 
\cite{S-loc,S-C2a,S-W2p}, where he recently introduced new techniques to obtain
global versions of all the previous regularity results 
under suitable regularity assumptions on the boundary data.
Let us describe the main results.

Assume that $0\in\partial\O$,
that $\O\subset \R^n$ is a bounded open convex set satisfying 
\begin{equation}
\label{eq:Omega zero}
B_\rho(\rho \,e_n)\subset \Omega\subset B_{1/\rho}(0)\cap \{x_n> 0\}
\end{equation}
for some $\rho>0$, and that $u:\O\to \R$ satisfies 
\begin{equation}
\label{eq:MA l L}
\mu_u =f\,dx\qquad \text{in }\O
\end{equation}
for some $0<\l\leq f \leq 1/\lambda$.
Extend $u$ by letting it being equal to $+\infty$ in $\R^n\setminus \overline\O$,
and up to subtracting a linear function assume that $\ell(x)\equiv 0$ is the tangent plane to $u$ at $0$,
that is
\begin{equation}
\label{eq:tangent}
u \geq 0,\qquad u(0)=0,\qquad \text{and}\qquad u(x)\not\geq \e\,x_n \quad \forall\,\e>0.
\end{equation}
The main result in \cite{S-loc} shows that
if $u \approx |x|^2$ along $\p\O \cap \{x_n \leq \rho\}$, then the
sections $S_t(0):=\{x \in \O\,:\,u(x)<t\}$ are comparable to half-ellipsoids for $t$ small.
More precisely, 
the following holds:
\begin{theo}
\label{thm:loc bdry}
Let $\O\subset \R^n$ be a bounded open convex set satisfying \eqref{eq:Omega zero},
and $u:\O\to \R$ be  a convex function satisfying \eqref{eq:MA l L} for some $0<\l\leq f \leq 1/\lambda$. Assume that \eqref{eq:tangent} holds, and that
\begin{equation}
\label{eq:quadr}
\beta\,|x|^2 \leq u(x)\leq \frac{1}{\beta}\,|x|^2 \qquad \text{on }\p\O\cap \{x_n \leq \rho\}
\end{equation}
for some $\beta >0$. Then, for any $t>0$ small,
there exists an ellipsoid $\E_t$ of volume $t^{n/2}$ such that
$$
\Bigl(\frac{1}{K}\,\E_t\Bigr)\cap \overline \O\subset S_t(0)\subset \bigl(K\,\E_t\bigr)\cap\overline \O,
$$
where $K>1$ depends only on $n,\l,\rho$, and $\b$.
In addition, the ellipsoid $\E_t$ is comparable to a ball of radius $\sqrt{t}$, up to a possible translation
along the $x_n$-direction of size $|\log t|$. Specifically, there exists a linear transformation $A_t:\R^n\to\R^n$
of the form
$$
A_t(x)=x-\tau\,x_n,\qquad \tau=(\tau_1,\ldots,\tau_{n-1},0)\in \R^n,\qquad\text{and}\qquad |\tau|\leq K\,|\log t|,
$$
such that $\E_t=A_t\bigl(B_{\sqrt{t}}(0)\bigr)$.
\end{theo}

The last part of the above result provides information about the behavior of the second derivatives near the origin. Indeed,
heuristically, this result states that inside $S_t(0)$ the tangential second derivatives are uniformly bounded 
both from above and below, while the mixed second derivatives are bounded by $|\log t|$.
This is very interesting given that $\mu_u$ is only bounded from above and below,
and that the boundary data as well as the boundary are only $C^{1,1}$.


As a consequence of Theorem~\ref{thm:loc bdry} and the interior estimates proved in Section \ref{sect:Caff}, in \cite{S-C2a,S-W2p} Savin obtained the following global $C^{2,\a}$-$W^{2,p}$ estimates.
\begin{theo}
Let $\O\subset \R^n$ be a bounded open uniformly convex set,
 $u:\O\to \R$ be a convex function satisfying \eqref{eq:MA l L} for some $0<\l\leq f \leq 1/\l$,
 and assume that both $u|_{\p\O}$ and $\p\O$ are of class $C^{1,1}$.
Suppose also that $u$ separates quadratically on $\p\O$ from its tangent plane, that is
$$
u(z)-u(x) \geq \<\n u(x),z-x\>+\beta\,|z-x|^2\qquad \forall\,x,z\in\p\O.
$$
Then:
\begin{enumerate}
\item[-] There exist $\e>0$ such that $u \in W^{2,1+\e}(\overline\O)$.
\item[-] For any $p>1$, if $\|f-1\|_{L^\infty(\overline\O)}\leq  e^{-Cp}$ then $u \in W^{2,p}(\overline\O)$.
\item[-] Assume that $f \in C^{0,\a}(\overline\O)$ and that 
both $u|_{\p\O}$ and $\p\O$ are of class $C^{2,\a}$. Then $u \in C^{2,\a}(\overline\O)$.
\end{enumerate}
\end{theo}

As observed in \cite{S-loc}, the assumption that
$u$ separates quadratically on $\p\O$ from its tangent plane
is verified, for instance, whenever $\p\O$ and $u|_{\p\O}$
are of class $C^3$ with $\O$ uniformly convex.

\section{Recent developments III: Smoothness of the first eigenfunction}

Let $\Omega$ be a smooth uniformly convex set.
In the paper \cite{Lions}, P.-L. Lions investigated the existence and uniqueness of the first eigenvalue for the Monge-Amp\`ere operator, namely the existence of a nontrivial convex function $\psi_1\in C^{1,1}(\overline\Omega)\cap C^\infty(\Omega)$ and a positive constant $\lambda_1$ 
such that
\begin{equation}
\label{eq:MA egn}
(\det D^2\psi_1)^{1/n}=-\lambda_1\psi_1\quad
\text{in }\Omega,\qquad \psi_1=0\quad
\text{on }\partial\Omega.
\end{equation}
As shown in \cite{Lions}, the couple $(\lambda_1,\psi_1)$ is essentially unique. More precisely, if $\psi:\Omega\to \R$ is a nontrivial convex
function and $\lambda$ a positive constant such that
$$
(\det D^2\psi)^{1/n}=-\l\psi\quad
\text{in }\Omega,\qquad \psi=0\quad
\text{on }\partial\Omega,
$$
then $\lambda=\lambda_1$ and $\psi=\theta\psi_1$ for some positive constant $\theta$.

Using the algebraic formula
$$
(\det A)^{1/n}=\inf\Bigl\{{\rm tr}(AB)\,:\,B\text{ symm. pos. def.,}\,\det B \geq \frac{1}{n^n} \Bigr\},
$$
one can prove that
$$
\l_1=\inf\Bigl\{\l_1(a_{ij})\,:\,a_{ij}\in C(\overline\Omega),\,a_{ij}\text{ symm. pos. def.,} \,\det(a_{ij})\geq \frac{1}{n^n} \Bigr\},
$$
where $\l_1(a_{ij})$ is the first eigenvalue of the linear elliptic operator $a_{ij}\partial_{ij}$.
In addition, thanks to this formula one can approximate the Monge-Amp\`ere equation with Hamilton-Jacobi-Bellman equations of the form
$$
\mathcal A^\epsilon\psi:=\inf\Bigl\{a_{ij}\p_{ij}\psi\,:\,a_{ij}\in C(\overline\Omega),\,a_{ij}\text{ symm. pos. def.,} \,\det(a_{ij})\geq \frac{1}{n^n},\,{\rm tr}(a_{ij})\leq \frac1{\epsilon}  \Bigr\},
$$
and deduce some interesting stochastic interpretation for $\lambda_1$
(see \cite{Lions2,Lions} for more details).

As observed in Section \ref{sect:Alex sol}, many results for the equation $\mu_u=f(x)\,dx$ can be extended to the general case $\mu_u=f(x,u,\nabla u)\,dx$ provided $\partial_uf \geq0$, as this ensures the validity of the maximum principle. An interesting consequence of Lions' result is the validity of a maximum principle also when $f$ is slightly decreasing with respect to $u$.
More precisely, the equation $\mu_u=F(x,u)\,dx$ has a unique solution provided $\partial_u \bigl(F(x,u)^{1/n}\bigr) > -\lambda_1$ (see \cite[Corollary 2]{Lions}).

\smallskip

Note that, in view of the $C^{1,1}$ regularity of $\psi_1$, near the boundary of $\Omega$ one can write $|\psi_1(x)|=g(x)d_{\partial\Omega}(x)$,
where $g:\Omega\to \R$ is a strictly positive Lipschitz function, and $d_{\partial\Omega}(x)={\rm dist}(x,\partial\Omega)$ denotes the distance function to the boundary.
In other words, $\psi_1$ solves a Monge-Amp\`ere equation of the form
\begin{equation}
\label{eq:MA deg}
\det D^2\psi_1=G\,d_{\partial\Omega}^n\quad
\text{in }\Omega,\qquad \psi_1=0\quad
\text{on }\partial\Omega,
\end{equation}
where $G\geq c_0>0$ is Lipschitz.

Because the right hand side vanishes on $\partial\Omega$, this equation is degenerate near the boundary and it has been an open problem for more than 30 years whether $\psi_1$ is smooth up to the boundary. The solution to this question has been given only recently, first by Hong, Huang, and Wang in two dimensions \cite{HHW}, and then by Savin \cite{S-loc2} and by Le and Savin \cite{LeSavin} in arbitrary dimensions.

More precisely, consider the general calss of Monge-Amp\`ere equations
\begin{equation}
\label{eq:MA deg2}
\mu_u=f(x)\,dx\quad
\text{in }\Omega,\qquad u=0\quad
\text{on }\partial\Omega,\qquad f(x)=G(x)d_{\partial\Omega}(x)^s,
\end{equation}
where  $s>0$ and $G$ is a continuous strictly positive function.
In \cite{S-loc2} Savin proved the following $C^2$ regularity estimate at the boundary:

\begin{theo}
\label{thm:loc bdry2}
Let $\O\subset \R^n$ be a bounded open convex set satisfying \eqref{eq:Omega zero},
and $u:\O\to \R$ be  a convex function satisfying \eqref{eq:MA deg2}. Assume that \eqref{eq:tangent} and \eqref{eq:quadr} hold, and that $u|_{\partial\Omega\cap B_\rho(0)}$ is of class $C^2$ for some $\rho>0$.
 Then $u$ is $C^2$ at $0$. More precisely, there exists a vector $\tau$ perpendicular to $e_n$, a quadratic polynomial $Q:\R^{n-1}\to \R$, and a constant $a>0$, such that
 $$
 u(x+\tau x_n)=Q_0(x')+ax_n^{2+s}+o(|x'|^2+x_n^{2+s})\qquad \forall\,x=(x',x_n) \in B_\rho(0).
 $$ 
\end{theo}
As a consequence of this result, since \eqref{eq:MA egn} is of the form \eqref{eq:MA deg2} with $s=n$, Savin obtained the global $C^2$ regularity of the first eigenfunction:
\begin{coro}
Let $\O\subset \R^n$ be a uniformly convex set of class $C^2$,
and let $\psi_1$ be the first eigenfunction (see \eqref{eq:MA egn}).
Then $\psi_1\in C^{2}(\overline\Omega).$
\end{coro}

By a perturbative approach based on Theorem \ref{thm:loc bdry2}, Le and Savin improved the boundary $C^2$ regularity to $C^{2,\beta}$.
More precisely, they showed the following pointwise estimate:

\begin{theo}
\label{thm:loc bdry3}
Let $\O\subset \R^n$ be a bounded open convex set satisfying \eqref{eq:Omega zero},
and $u:\O\to \R$ be  a convex function satisfying \eqref{eq:MA deg2}. Assume that \eqref{eq:tangent} and \eqref{eq:quadr} hold, and that $u|_{\partial\Omega\cap B_\rho(0)}$ is of class $C^{2,\beta}$ for some $\beta \in \left(0,\frac{2}{2+s}\right)$ and $\rho>0$.
 Also, assume that $G \in C^{0,\gamma}(\overline\Omega\cap B_\rho(0))$ for some $\gamma \geq \frac{\beta(2+s)}{2}$.
 Then $u$ is $C^{2,\beta}$ at $0$. More precisely, there exists a vector $\tau$ perpendicular to $e_n$, a quadratic polynomial $Q:\R^{n-1}\to \R$, and a constant $a>0$, such that
 $$
 u(x+\tau x_n)=Q_0(x')+ax_n^{2+s}+O(|x'|^2+x_n^{2+s})^{1+\beta/2}\qquad \forall\,x=(x',x_n) \in B_\rho(0).
 $$ 
\end{theo}
As a consequence of this result, one obtains the 
global $C^{2,\beta}$ regularity of the first eigenfunction for any $\beta<\frac{2}{2+n}$.

We note that usually, in this type of elliptic questions, once one obtains  $C^{2,\beta}$ regularity then the higher regularity follows easily by Schauder estimates. This is not the case in this situation because of the high degeneracy of the equation. The key idea in \cite{LeSavin} consists in performing both an hodograph transform and a partial Legendre transform in order to deduce that (a suitable transformation of) the first eigenfunction satisfies a degenerate Grushin-type equation with H\"older coefficients.
Once this is achieved, Le and Savin conclude the global smoothness of $\psi_1$ by applying
Schauder estimates for Grushin-type operators:
\begin{coro}
Let $\O\subset \R^n$ be a uniformly convex set of class $C^\infty$,
and let $\psi_1$ be the first eigenfunction (see \eqref{eq:MA egn}).
Then $\psi_1\in C^{\infty}(\overline\Omega).$
\end{coro}

\bigskip

{\it Acknowledgments:} The author is thankful to Yash Jhaveri and Connor Mooney for useful comments on a preliminary version of the manuscript.
The author is supported by the ERC Grant ``Regularity and Stability in Partial Differential Equations (RSPDE)''.

\end{document}